\documentclass[pdflatex, sn-mathphys-num, iicol]{sn-jnl}

\usepackage{graphicx}
\usepackage{multirow}
\usepackage{amsmath,amssymb,amsfonts}
\usepackage{amsthm}
\usepackage{mathrsfs}
\usepackage[title]{appendix}
\usepackage{xcolor}
\usepackage{textcomp}
\usepackage{manyfoot}
\usepackage{booktabs}
\usepackage{algorithm}
\usepackage{algorithmicx}
\usepackage{algpseudocode}
\usepackage{listings}
\usepackage{enumitem}

\theoremstyle{thmstylethree}

\newtheorem{definition}{Definition}
\numberwithin{definition}{section}

\numberwithin{remark}{section}

\begin{document}

\title
{A framework for generating nonconforming triangular meshes with multiple discretization layers}

\author{\fnm{I.\,L.}\sur{Semenov}}
\author{\fnm{M.\,M.}\sur{Becker}}

\affil{
\orgname{Leibniz Institute for Plasma Science and Technology},
\orgaddress{
\street{Felix-Hausdorff-Str.~2},
\city{Greifswald}, 
\postcode{17489},
\country{Germany}}}

\abstract{
We present a framework for generating nonconforming triangular meshes with multiple discretization layers. The framework exploits characteristic structural properties of meshes produced by a frontal Delaunay algorithm with uniform element size. The bulk region of such meshes exhibits a structured pattern resembling a regular triangular lattice. Owing to this structure, the bulk region can be coarsened by grouping elements into connected subsets of larger composite elements. When applied repeatedly, this procedure produces a mesh with multiple discretization layers. For complex geometries, the framework can be used to create composite multidomain meshes, where multiple discretization layers are generated within each subdomain.
We present examples of meshes generated by the proposed framework and discuss postprocessing strategies for their refinement and coarsening. The resulting meshes are well suited for the application of adaptive mesh refinement techniques. The proposed framework can be readily integrated with existing mesh generators and finite-element solvers that support nonconforming triangular meshes.}

\keywords{nonconforming triangular meshes, mesh generation, frontal Delaunay algorithm, adaptive mesh refinement}
\maketitle
\section{Introduction}

Nonconforming (irregular) meshes are meshes in which elements do not necessarily share complete common edges or faces. Such meshes are typically used in the context of adaptive mesh refinement (AMR) and are supported by most modern finite-element solvers \citep{cerveny2019, anderson2021mfem, dobrev2022, arndt2023deal}.

The most common form of nonconforming AMR is implemented on logically structured meshes using hierarchical data structures such as quadtrees or octrees~\citep{burstedde2011p4est,teunissen2018afivo}.
On unstructured meshes, nonconforming AMR can be realized using the red-refinement approach~\citep{carey1997, giuliani2019adaptive, solin2023selected}, in which elements of the background mesh are subdivided into subelements that can later be coarsened. This approach, however, is less efficient than in the structured case when the background mesh does not possess a suitable hierarchy of discretization layers. Creating such layers therefore requires the development of new algorithmic approaches that operate in fully unstructured settings.

Nonconforming meshes also provide a flexible alternative to nonuniform conforming meshes for resolving geometrically complex regions of the computational domain. In this regard, a transition from logically structured to unstructured settings for nonconforming mesh generation offers greater flexibility as well.

In an unstructured setting, generating a nonconforming mesh with multiple discretization layers involves two key aspects. First, it is necessary to identify an appropriate data structure capable of representing nonconforming element interfaces within a mesh. Such a data structure should be sufficiently flexible to allow the straightforward implementation of basic combinatorial operations on mesh elements. Second, it is necessary to provide an approach for constructing a hierarchy of discretization layers within the mesh that resolves the geometry of the computational domain while aiming to minimize the total number of mesh elements. Such a discretization can then serve as a basis for AMR techniques. 

In this paper, we address the aforementioned aspects in the context of triangular meshes. We present a framework for generating nonconforming triangular meshes with multiple discretization layers in domains with relatively simple geometry. The proposed framework builds on the following observation. The application of a frontal Delaunay algorithm~\citep{remacle2013frontal} to triangular mesh generation with a uniform element size results in meshes exhibiting a characteristic two-region structure. The first region is an unstructured layer along the domain boundary that resolves the geometric complexity of the domain. This layer is connected to a bulk region that closely resembles a regular triangular lattice. Due to its regular structure, the bulk region can be coarsened by grouping elements into connected subsets of larger composite elements. When applied repeatedly, this procedure produces a nonconforming triangular mesh with multiple discretization layers.

The characteristic structure of the background mesh described above is generally observed for domains with simple geometry. For more complex geometries, the framework can be used to create composite multidomain meshes, where multiple discretization layers are generated within each subdomain.

The practical implementation of the proposed framework relies on representing nonconforming element interfaces as zero-area triangles. With this approach, the mesh remains conforming in the combinatorial sense. From a geometric perspective, a nonconforming mesh is treated as a limiting case of a conforming mesh, where nonconforming element interfaces correspond to triangles of infinitesimally small area. This approach unifies the description of conforming and nonconforming triangular meshes and enables straightforward implementation of basic combinatorial operations on mesh elements.

We demonstrate the application of the proposed framework to the generation of meshes with multiple discretization layers in simple geometric domains and discuss an approach for subsequent refinement and coarsening of the resulting meshes.
We also demonstrate how the proposed framework can be used to construct composite multidomain meshes with multiple discretization layers.

The structure of the paper is as follows. In Section~\ref{sec:MeshDef} we introduce the representation of nonconforming meshes adopted in this work.
In Section~\ref{sec:SupTri} we introduce a specific data structure, called a {\it supertriangulation}, which is used in the mesh generation algorithm. The mesh generation algorithm is presented in Section~\ref{sec:MeshConstr}. Mesh postprocessing procedures, such as refinement and coarsening, are discussed in Section~\ref{sec:MeshProcess}. Examples of mesh generation are presented in Section~\ref{sec:Examps}. The concluding remarks are given in Section~\ref{sec:Concls}.

\section{Mesh representation}
\label{sec:MeshDef}

In the present work, nonconforming triangular meshes are treated as a limiting case of conforming meshes, where nonconforming interfaces are represented by triangles of infinitesimally small area. The idea behind this approach is discussed below.

\subsection{Preliminary considerations}

\begin{figure*}[!t]
\centering
\includegraphics[scale=1]{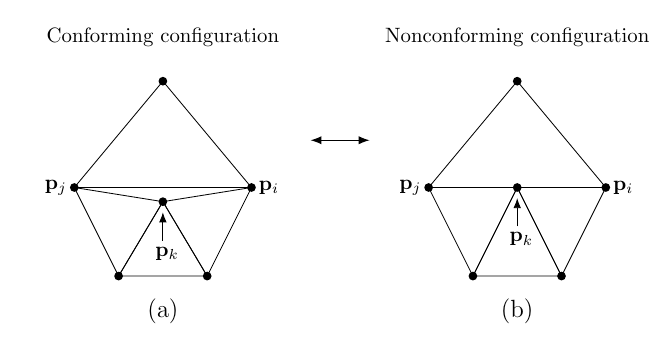}
\caption{\label{fig:TriuExamp}
Transition between conforming and nonconforming mesh configurations. Panel (a) shows a conforming configuration in which a nonconforming interface is represented as a triangle of infinitesimally small area. Panel (b) illustrates the limiting case of the configuration in panel (a). In this case, the nonconforming interface is represented by a triangle of zero area.}
\end{figure*}

Let \(\mathcal{P} = \{\mathbf{p}_i\}_{i=1}^n \subset \mathbb{R}^2\) be a finite set of \(n \ge 3\) distinct points in the plane, and let \(\Omega = \mathrm{conv}(\mathcal{P})\) denote its convex hull. Assume that \(\Omega\) has nonzero area and admits a conforming triangulation,
\begin{equation}
\label{eq:Triu}
\mathcal{T} = 
\left\{
\triangle(\mathbf{p}_{i_\alpha}, \mathbf{p}_{j_\alpha}, \mathbf{p}_{k_\alpha})
\right\}_{\alpha = 1}^m,
\end{equation}
where \(\mathbf{p}_{i_\alpha}, \mathbf{p}_{j_\alpha}, \mathbf{p}_{k_\alpha} \in \mathcal{P}\), as defined, for example, in~\citep{frey2008mesh}. The vertices of each triangle in Eq.~\eqref{eq:Triu} are listed in \emph{counterclockwise} order, and the operator \(\triangle(\cdot, \cdot, \cdot)\) denotes the convex hull of its arguments.

Suppose now that \(\mathcal{T}\) contains a nearly degenerate triangle with vertices \(\mathbf{p}_i\), \(\mathbf{p}_j\), and \(\mathbf{p}_k\), as illustrated in Fig.~\ref{fig:TriuExamp}(a). Geometrically, \(\mathbf{p}_k\) can be continuously moved toward the configuration shown in Fig.~\ref{fig:TriuExamp}(b), while preserving the connectivity of the triangulation. In this limiting configuration, \(\mathbf{p}_k\) becomes collinear with \(\mathbf{p}_i\) and \(\mathbf{p}_j\), forming a typical nonconforming interface in which \(\mathbf{p}_k\) serves as a hanging node.

This observation suggests that nonconforming meshes can be interpreted as a limiting case of conforming meshes, where nonconforming interfaces are represented by triangles of infinitesimally small area. A comprehensive theoretical formulation of this transition is beyond the scope of the present work.

Nevertheless, the practical implementation of this idea is straightforward.

\subsection{Unified mesh description}

In practice, triangulations are implemented using data structures that represent triangles by referencing the indices of vertices in the underlying point set. A standard example of such a representation is the \emph{connectivity matrix}. For the triangulation \(\mathcal{T}\) in Eq.~\eqref{eq:Triu}, this matrix is
\begin{equation*}
\mathsf{T}_{\mathcal{T}} = 
\begin{bmatrix}
i_1 & j_1 & k_1 \\
\vdots & \vdots & \vdots \\
i_m & j_m & k_m \\
\end{bmatrix} \in \mathbb{N}^{m \times 3},
\end{equation*}
where each row contains the indices of the points in \(\mathcal{P}\) that form a triangle in \(\mathcal{T}\).


In this work, the connectivity matrix is used to provide a unified representation of both conforming and nonconforming triangular meshes. Nonconforming interfaces are incorporated into the matrix as additional rows (elements), using a vertex ordering that preserves mesh conformity at the combinatorial level. 

For example, the nonconforming interface shown in Fig.~\ref{fig:TriuExamp}(b) appears in the connectivity matrix as a standard triangle:
\[
\mathsf{T}_{\mathcal{T}} = 
\begin{bmatrix}
\vdots & \vdots & \vdots \\
i & j & k \\
\vdots & \vdots & \vdots
\end{bmatrix}.
\]
At the combinatorial level, this representation corresponds to the conforming geometric configuration depicted in Fig.~\ref{fig:TriuExamp}(a) and captures the limiting transition to the nonconforming arrangement depicted in Fig.~\ref{fig:TriuExamp}(b). 

If we assume that the connectivity matrix remains unchanged during the continuous transition from configuration~(a) to configuration~(b), then the area \(A\) of the triangle \(\triangle(\mathbf{p}_i, \mathbf{p}_j, \mathbf{p}_k)\) tends to zero from above, that is, \(A \to 0^+\). Thus, the nonconforming configuration~(b) is approached as the limit of a continuous sequence of well-defined, infinitesimally close conforming configurations.

This picture underlies the basic idea of our approach. We use the connectivity matrix as the primary data structure for mesh representation. During mesh construction, nonconforming interfaces are added to the connectivity matrix in such a way that the mesh preserves conformity in the sense illustrated in Fig.~\ref{fig:TriuExamp}. That is, the mesh always admits a conforming triangulation that is infinitesimally close to a nonconforming one.

At the geometric level, the mesh is still represented by a collection of the form~\eqref{eq:Triu}, where nonconforming interfaces appear as degenerate triangles. We refer to these additional elements as \emph{void elements}, since they have zero area. For simplicity, we consider only a specific class of void elements in which
\begin{equation}
\label{eq:MidNodeVoid}
\mathbf{p}_k = \tfrac{1}{2}( \mathbf{p}_i + \mathbf{p}_j).
\end{equation}
The node $\mathbf{p}_k$ is referred to as the \emph{hanging node}.

Note that in our work, all operations on the mesh are performed at the combinatorial level using the connectivity matrix, while the geometric representation is used solely for visualization.

\section{Supertriangulation}
\label{sec:SupTri}

Operations involving nonconforming meshes are commonly encountered in the context of adaptive mesh refinement, which inherently involves multiple discretization levels. Supporting such levels on triangular meshes requires additional data structures. In particular, grouping elements into top-level complexes is crucial for hierarchical mesh management. To this end, we introduce a data structure called the \emph{supertriangulation}.

\subsection{Basic definitions}
\label{sec:SupTriDef}

\begin{figure*}[!h]
\centering
\includegraphics[width=0.8\linewidth]{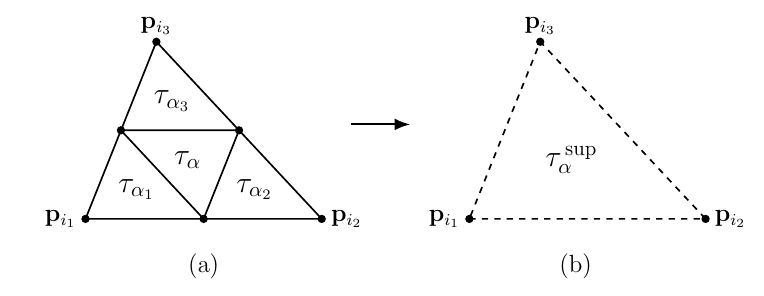}
\caption{\label{fig:SupTriu}
Illustration of the supertriangle definition. Panel (a) shows a group of elements in the background mesh used to define a supertriangle. Panel (b) shows the supertriangle formed from the group of elements shown in panel (a).}
\end{figure*}

\begin{figure*}[!t]
\centering
\includegraphics[width=0.8\linewidth]{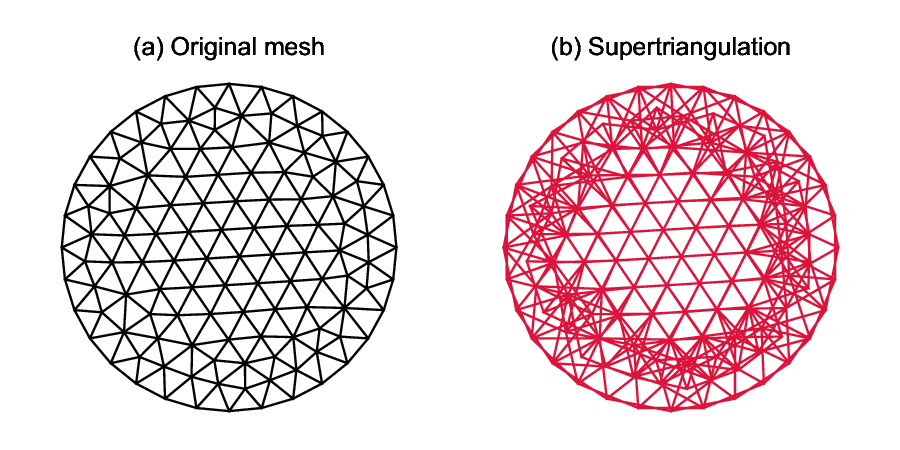}
\caption{\label{fig:SupTriExamp}
Illustration of a supertriangulation of a mesh generated in the unit circle. Panel (a) shows the original background mesh and panel (b) shows the supertriangulation constructed from this mesh.}
\end{figure*}

Let \(\mathcal{T}\) be a triangulation of the point set \(\mathcal{P}\), as introduced in Section~\ref{sec:MeshDef}. We define \(\mathcal{T}_{\mathrm{core}}\) to be the subset of triangles in \(\mathcal{T}\) whose edges are all shared with adjacent triangles in the mesh. For later reference, we introduce the operator
\begin{equation}
\label{eq:Tcore}
T^{\mathrm{core}}(\mathcal{T}) := \mathcal{T}_{\mathrm{core}}.
\end{equation}

Consider a triangle \(\tau_\alpha \in \mathcal{T}_{\text{core}}\), as illustrated in Fig.~\ref{fig:SupTriu}(a). Its neighboring triangles are denoted by \(\tau_{\alpha_1}, \tau_{\alpha_2}, \tau_{\alpha_3}\), listed in counterclockwise order. Let \(\mathbf{p}_{i_k}\) denote the vertex of \(\tau_{\alpha_k}\) that lies opposite the edge it shares with \(\tau_\alpha\), for each \(k = 1, 2, 3\).

We define the corresponding \emph{supertriangle} associated with \(\tau_\alpha\) as
\begin{equation}
\label{eq:SupTriang}
\tau_\alpha^{\mathrm{sup}} := \triangle(\mathbf{p}_{i_1}, \mathbf{p}_{i_2}, \mathbf{p}_{i_3}),
\end{equation}
as illustrated in Fig.~\ref{fig:SupTriu}(b).

Building on this, we introduce the operator \(T^{\mathrm{sup}}\), which maps a subset of \(\mathcal{T}_{\mathrm{core}}\) to its corresponding set of supertriangles:
\begin{equation}
\label{eq:Tsup}
T^{\mathrm{sup}}(\mathcal{T}_{\mathrm{core}}') := 
\left\{ \tau_\alpha^{\mathrm{sup}} \;\middle|\; \tau_\alpha \in \mathcal{T}_{\mathrm{core}}' \right\}
\end{equation}
where \(\mathcal{T}_{\mathrm{core}}' \subseteq \mathcal{T}_{\mathrm{core}}\).

The output of the operator \(T^{\mathrm{sup}}\) is referred to as a \emph{supertriangulation}. The set of all supertriangulations, denoted by \(\mathscr{T}_{\mathrm{sup}}\), is thus given by the image of the operator:
\begin{equation*}
\mathscr{T}_{\mathrm{sup}} := \mathrm{Im}\left( T^{\mathrm{sup}} \right).
\end{equation*}

The supertriangulation generated by applying \(T^{\mathrm{sup}}\) to \(\mathcal{T}_{\mathrm{core}}\) is referred to as the \emph{background supertriangulation}.

Figure~\ref{fig:SupTriExamp} illustrates a background supertriangulation for a mesh generated within the unit circle using the frontal Delaunay algorithm~\citep{remacle2013frontal} in Gmsh~\citep{geuzaine2009gmsh}.
It can be seen that supertriangles intersect geometrically, with irregular overlaps near the boundary of the mesh and more structured intersections in the interior region.

By construction, intersections between supertriangles can be interpreted purely in combinatorial terms. We therefore introduce the following definitions.

\begin{definition}
\label{def:CombStenc}
The \emph{combinatorial stencil} of a supertriangle is the set of indices of original mesh elements used to construct it.
For instance, the combinatorial stencil of the supertriangle \(\tau_\alpha^{\mathrm{sup}}\) in Fig.~\ref{fig:SupTriu} is given by \(\{\alpha, \alpha_1, \alpha_2, \alpha_3\}\).
\end{definition}

\begin{definition}
\label{def:SupIntersect}
Two supertriangles are said to be \emph{combinatorially disjoint} if their combinatorial stencils have no indices in common. Otherwise, they are said to be \emph{combinatorially overlapping}.
\end{definition}

\begin{definition}
\label{def:SupCompact}
A supertriangulation is said to be \emph{compact} if all of its supertriangles are pairwise combinatorially disjoint.
\end{definition}

\subsection{Supported operations}
\subsubsection{Filtering procedures}

According to Eq.~\eqref{eq:SupTriang}, a supertriangulation consists of triangle elements. Their connectivity can be characterized using a connectivity matrix, defined analogously to the connectivity matrix of the original triangulation.

As a result, supertriangles can be classified according to the number of their connections with other supertriangles. For later use, we introduce several procedures to filter supertriangles based on their connectivity.

\begin{definition}
We define the \emph{smoothing} of a supertriangulation \(\mathcal{T}_{\mathrm{sup}} \in \mathscr{T}_{\mathrm{sup}}\) as the process of removing all supertriangles with fewer than two connections to other supertriangles. This process can be performed in a single pass or iteratively, until no further removals are possible.
\end{definition}

\begin{definition}
We define the \emph{stripping} of a supertriangulation \(\mathcal{T}_{\mathrm{sup}} \in \mathscr{T}_{\mathrm{sup}}\) as the procedure that removes all supertriangles connected to exactly two other supertriangles. This operation is always performed in a single pass.
\end{definition}

\begin{definition}
\label{def:SupClean}
The \emph{cleaning} of a supertriangulation \(\mathcal{T}_{\mathrm{sup}} \in \mathscr{T}_{\mathrm{sup}}\) is defined as the process consisting of a stripping step followed by iterative smoothing. We denote this process by the operator
\[
F^{\mathrm{clean}} : \mathscr{T}_{\mathrm{sup}} \to \mathscr{T}_{\mathrm{sup}}.
\]
\end{definition}

The stripping and smoothing procedures eliminate peripheral supertriangles.
In some cases, it is also beneficial to detach the supertriangulation from the boundary of the original triangulation. We therefore introduce the \emph{detachment} procedure.

\begin{definition}
\label{def:Detach}
A supertriangulation \(\mathcal{T}_{\mathrm{sup}}\) associated with a triangulation \(\mathcal{T}\) is said to be \emph{detached} if all supertriangles with vertices lying on the boundary of \(\mathcal{T}\) have been removed. We denote this process by the operator
\[
F^{\mathrm{detach}} : \mathscr{T}_{\mathrm{sup}} \to \mathscr{T}_{\mathrm{sup}},
\]
which maps a supertriangulation to its detached version
\end{definition}

\begin{figure*}[!t]
\centering
\includegraphics[width=0.8\linewidth]{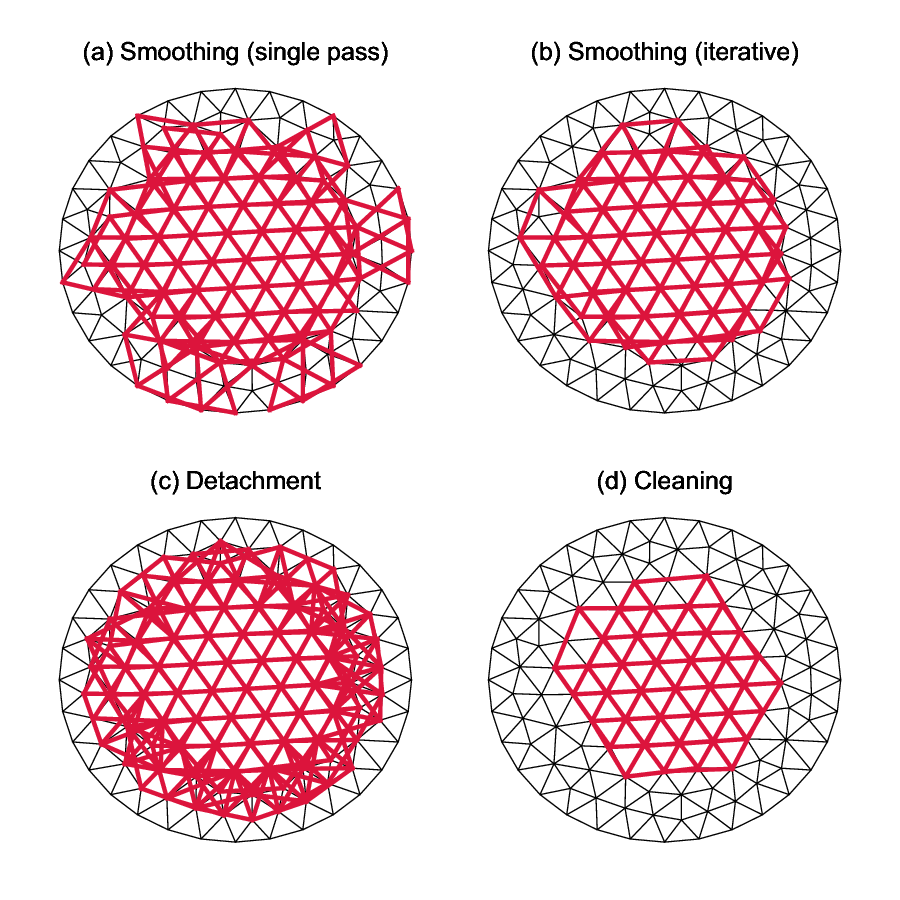}
\caption{\label{fig:SupTriFilt}
Illustration of filtering procedures applied to the supertriangulation of a mesh generated in the unit circle. The black color indicates the elements of the original background mesh. The red color indicates the elements of the supertriangulation filtered using single-pass smoothing (a), iterative smoothing (b), boundary detachment (c), and cleaning (d).
}
\end{figure*}

Figure~\ref{fig:SupTriFilt} illustrates all the filtering procedures defined above, as applied to the supertriangulation shown in Fig.~\ref{fig:SupTriExamp}.

\subsubsection{Reduction procedure}

An essential operation on a supertriangulation is the reduction to a connected subset of supertriangles. This process is carried out in two steps:

\begin{enumerate}[label=\arabic*)]
\item Construct the connectivity graph between all supertriangles in the supertriangulation.
\item Extract connected supertriangles via standard graph traversal.
\end{enumerate}

The connectivity graph is encoded by the binary adjacency matrix \(\mathsf{A}^{\mathrm{sup}}\), where each node represents a supertriangle. The element \(\mathsf{A}^{\mathrm{sup}}_{i,j}\) of the adjacency matrix is 1 if the \(i\)-th supertriangle is connected to the \(j\)-th supertriangle, and 0 otherwise.

Our implementation performs a breadth-first search on the adjacency matrix, starting from a preselected row. The reduction procedure is summarized in the following definition.

\begin{definition}
\label{def:ReducOpr}
The \emph{reduction operator} is defined as the mapping
\[
F^{\mathrm{reduce}} : \mathscr{T}_{\mathrm{sup}} \to \mathscr{T}_{\mathrm{sup}},
\]
which takes a supertriangulation \(\mathcal{T}_{\mathrm{sup}}\) and returns its filtered version containing only the supertriangles reachable by traversing the connectivity graph of \(\mathcal{T}_{\mathrm{sup}}\) from a preselected supertriangle \(\tau_s^{\mathrm{sup}} \in \mathcal{T}_{\mathrm{sup}}\).
\end{definition}

Figure~\ref{fig:SupTri_Compr} illustrates the results of the reduction procedure applied to the reference setup shown in Fig.~\ref{fig:SupTriExamp}. Two configurations are considered: one where the initial supertriangulation is detached from the mesh boundary, and one where it remains connected. The reduction process begins from the supertriangle whose centroid lies closest to the mean position of all supertriangle centroids.

\subsubsection{Iterative reduction}
\label{sec:IterReduc}

\begin{figure*}[!t]
\centering
\includegraphics[width=0.8\linewidth]{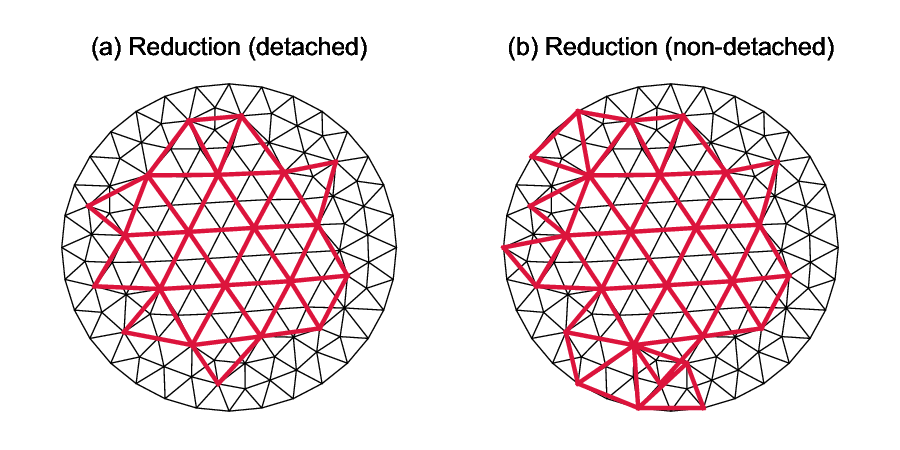}
\caption{\label{fig:SupTri_Compr}
Illustration of the reduction procedure applied to the supertriangulation shown in Fig.~\ref{fig:SupTriu}. Black color indicates the elements of the original background mesh and the red color indicates the elements of the reduced supertriangulation. Panel (a) shows the case in which the supertriangulation is detached from the mesh boundary before reduction, while panel (b) shows the supertriangulation after reduction without boundary detachment.
}
\end{figure*}

The output of the reduction process is a supertriangulation consisting of connected elements. However, combinatorial overlaps (see Definition~\ref{def:SupIntersect}) between elements are still possible, as demonstrated, for example, in Fig.~\ref{fig:SupTri_Compr}(b). To resolve this, we combine reduction with cleaning in an iterative procedure.

The algorithm iteratively applies the reduction step to the supertriangulation. If the result is not a compact supertriangulation, a cleaning step is executed, followed by another reduction attempt. This process continues until the algorithm produces either a compact or an empty supertriangulation. 

The complete algorithmic definition of this procedure is given below:

\begin{algorithm}[H]
\caption{Iterative Reduction}
\label{alg:IterReduc}
\begin{algorithmic}[1]
\State \textbf{Input:} Supertriangulation \(\mathcal{T}_{\mathrm{sup}}\)
\While{true}
    \State \(\mathcal{T}_{\mathrm{sup}} \gets F^{\mathrm{reduce}}(\mathcal{T}_{\mathrm{sup}})\)
    \If{\(\mathcal{T}_{\mathrm{sup}}\) is compact}
        \State \textbf{return} \(\mathcal{T}_{\mathrm{sup}}\)
    \Else
        \State \(\mathcal{T}_{\mathrm{sup}} \gets F^{\mathrm{clean}}(\mathcal{T}_{\mathrm{sup}})\)
        \If{\(\mathcal{T}_{\mathrm{sup}} = \emptyset\)}
            \State \textbf{return} \(\emptyset\)
        \EndIf
    \EndIf
\EndWhile
\end{algorithmic}
\end{algorithm}

To refer to the iterative reduction procedure, we introduce the following operator.

\begin{definition}
\label{def:IterReduc}
We define the \emph{iterative reduction operator} as the mapping
\[
F^{\mathrm{reduce}}_{\mathrm{iter}} : 
\mathscr{T}_{\mathrm{sup}} \to \mathscr{T}_{\mathrm{sup}},
\]
which takes a supertriangulation and applies the iterative reduction procedure as described in Algorithm~\ref{alg:IterReduc}.
\end{definition}

\section{Mesh construction}
\label{sec:MeshConstr}

The framework for supertriangulations introduced in Section~\ref{sec:SupTri} supports the construction of meshes with multiple discretization layers by enabling the extraction of coarse submeshes from the original mesh.

The unified mesh representation introduced in Section~\ref{sec:MeshDef} facilitates the construction of nonconforming triangular meshes. By treating nonconforming interfaces as regular triangles, it allows for greater flexibility in mesh manipulation. 

As mentioned early, the only essential requirement is the consistent inclusion of void elements that represent nonconforming interfaces into the connectivity matrix, to preserve the overall mesh integrity.

Utilizing these two features, we introduce an algorithm that generates a nonconforming triangular mesh by iteratively coarsening a high-quality conforming mesh.

Conforming triangular meshes can be automatically generated using modern mesh generators such as Gmsh \citep{geuzaine2009gmsh}. In particular, high-quality, nearly uniform meshes can be produced using the frontal Delaunay algorithm \citep{remacle2013frontal}. These meshes typically exhibit a composite structure, consisting of an unstructured boundary layer and a structured bulk that resembles a regular triangular lattice. This bulk can then be iteratively coarsened to produce a nonconforming triangular mesh with multiple discretization layers. Alternatively, the coarsening procedure can be applied to a structured triangular mesh that supports a connected supertriangulation pattern.

To illustrate the steps of the algorithm, we again use the conforming mesh shown in Fig.~\ref{fig:SupTriExamp} as a representative example of the input configuration. It is important to note, however, that the algorithm can operate on both conforming and nonconforming meshes. This flexibility enables its iterative use to construct triangular meshes with multiple discretization levels.

The structure of this section is as follows. In Subsections~\ref{sec:InputConfig}, \ref{sec:MeshDecomp}, and \ref{sec:MeshReconn}, we present the main stages of the mesh-coarsening algorithm. Subsection~\ref{sec:AlgMeshSumm} provides a summary of the algorithm.

\subsection{Input processing}
\label{sec:InputConfig}

Suppose we are given a triangulation \(\mathcal{T}\) over a point set \(\mathcal{P}\), as introduced in Section~\ref{sec:MeshDef}. This triangulation serves as the input to the algorithm and may contain both conforming and nonconforming element interfaces. The input processing step aims to construct a compact supertriangulation associated with \(\mathcal{T}\).

\subsubsection{Background supertriangulation}
\label{sec:InputBkgSup}

For this, we first construct the background supertriangulation, 
\begin{equation}
\label{eq:T_bg}
\mathcal{T}^{\,\rm bgk}_{\rm sup} = 
T^{\rm sup}\left(T^{\rm core}\left(\mathcal{T}\right)\right),
\end{equation}
where the operators $T^{\rm core}$ and $T^{\rm sup}$ are defined in Eqs.~\eqref{eq:Tcore} and \eqref{eq:Tsup}, respectively.

The background supertriangulation may optionally be detached from the boundary of \(\mathcal{T}\) as
\begin{equation}
\label{eq:T_bg_d}
\widetilde{\mathcal{T}}^{\,\rm bgk}_{\rm sup} = 
F^{\rm detach}\left(\mathcal{T}^{\,\rm bgk}_{\rm sup}\right),
\end{equation} 
where the detachment operator is introduced in Definition~\ref{def:Detach}.

\subsubsection{Reduction step}
\label{sec:InputIterRed}

We then apply the iterative reduction operator introduced in Definition~\ref{def:IterReduc} to \(\mathcal{T}^{\,\rm bgk}_{\rm sup}\) or \(\widetilde{\mathcal{T}}^{\,\rm bgk}_{\rm sup}\): 
\begin{equation}
\mathcal{T}_{\mathrm{sup}}^{\,\rm reduced} = 
F^{\mathrm{reduce}}_{\mathrm{iter}} 
\left(
	\mathcal{T}_{\rm sup}^{\,\rm input}
\right),
\end{equation}
where
\begin{equation*}
\mathcal{T}_{\rm sup}^{\,\rm input} = 
\begin{cases}
\widetilde{\mathcal{T}}^{\,\rm bgk}_{\rm sup},
& \text{ if detached},
\\
\mathcal{T}^{\,\rm bgk}_{\rm sup},
& \text{ otherwise},
\end{cases}
\end{equation*}
with \(\mathcal{T}_{\rm sup}^{\,\rm bgk}\) and \(\widetilde{\mathcal{T}}_{\rm sup}^{\,\rm bgk}\) given by Eqs.~\eqref{eq:T_bg} and~\eqref{eq:T_bg_d}, respectively.

The starting supertriangle in the reduction process is determined by selecting the supertriangle whose centroid lies closest to the anchor point $\mathbf{c}_0$, wihch is provided as the control parameter of the algorithm.

We assume that the iterative reduction ultimately yields a non-empty, compact supertriangulation \(\mathcal{T}_{\mathrm{sup}}^{\,\rm reduced}\). Otherwise, the process is terminated.

\subsubsection{Cleaning step}
\label{sec:InputCleaning}

In addition, we assume that the cleaning procedure may be applied to \(\mathcal{T}_{\,\mathrm{sup}}^{\,\rm reduced}\) multiple times as a post-reduction step. For instance, Fig.~\ref{fig:MeshCoarse}(a) shows the reduced supertriangulation from Fig.~\ref{fig:SupTri_Compr}(a) after a single cleaning step. This option is introduced to control the separation and size of discretization layers within the mesh. We assume that the supertriangulation remains non-empty during the cleaning process, otherwise the mesh-coarsening procedure is terminated.

In summary, the input configuration of the mesh-coarsening procedure consists of the original triangulation \(\mathcal{T}\) and a non-empty, compact supertriangulation \(\mathcal{T}_{\rm sup}^{\,\rm reduced}\) derived from it.

\subsection{Mesh decomposition}
\label{sec:MeshDecomp}

\begin{figure*}[!t]
\centering
\includegraphics[width=\linewidth]{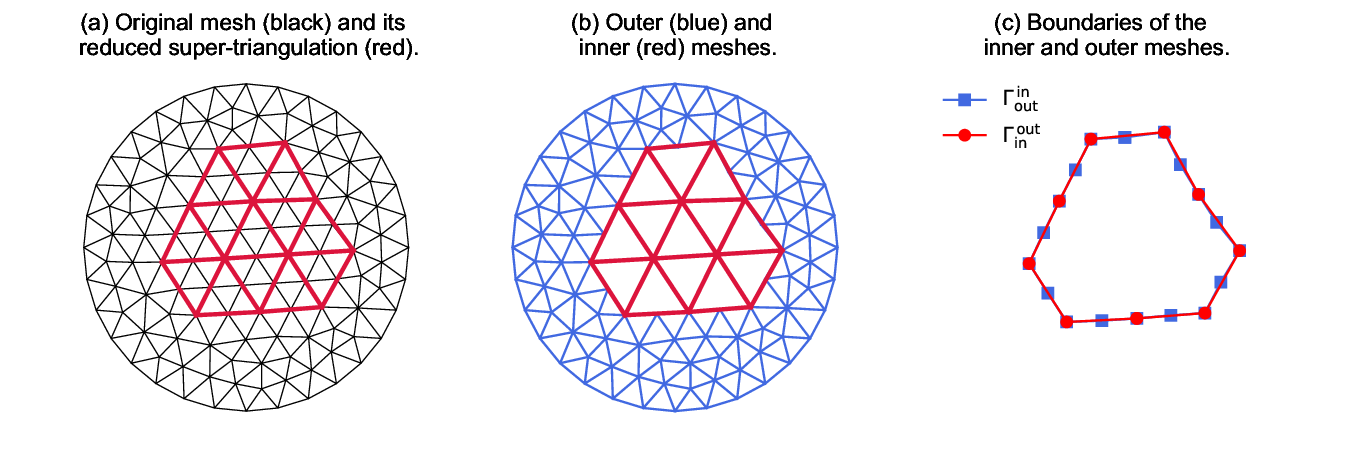}
\caption{\label{fig:MeshCoarse}
Illustration of the mesh decomposition step. Panel (a) shows the original background mesh together with its reduced supertriangulation. Panel (b) shows the outer and inner meshes constructed from the configuration shown in panel (a). Panel (c) shows the polygonal paths defining the contact between the inner and outer meshes.
}
\end{figure*}

Starting from the input configuration, we decompose the triangulation \(\mathcal{T}\) into two distinct components: the \emph{inner mesh} and the \emph{outer mesh}.

The inner mesh, which we denote by $\mathcal{T}_{\mathrm{in}}$, is the reduced supertriangulation $\mathcal{T}_{\rm sup}^{\,\rm reduced}$ itself.

The \emph{outer mesh}, which we denote by \(\mathcal{T}_{\mathrm{out}}\), is constructed by removing from \(\mathcal{T}\) all triangles covered by the reduced supertriangulation \(\mathcal{T}_{\rm sup}^{\,\rm reduced}\), namely:
\[
\mathcal{T}_{\mathrm{out}} = \mathcal{T} \setminus 
\widetilde{\mathcal{T}}_{\rm sup}^{\,\rm reduced},
\]
where
\begin{equation*}
\widetilde{\mathcal{T}}_{\rm sup}^{\,\rm reduced} = 
\{ \tau_\alpha \in \mathcal{T} : \alpha \in \mathcal{C}_{\rm sup}^{\,\rm reduced} \}
\end{equation*}
with \(\mathcal{C}_{\rm sup}^{\,\rm reduced}\) being the union of the combinatorial stencils of all supertriangles in \(\mathcal{T}_{\rm sup}^{\,\rm reduced}\) (see Definition~\ref{def:CombStenc}). Figure~\ref{fig:MeshCoarse}(b) shows the inner and outer meshes corresponding to the initial configuration in Fig.~\ref{fig:MeshCoarse}(a).

Algorithmically, we restrict attention to configurations in which \(\mathcal{T}_{\mathrm{out}}\) has an annular topology, that is, its boundary consists of a distinct inner and outer component. If this condition is not satisfied, the mesh-coarsening procedure is terminated.

We denote the inner boundary of the outer mesh by
\begin{equation}
\label{eq:GammaOutIn}
\Gamma_{\rm out}^{\rm in} := (\mathbf{q}_0, \ldots, \mathbf{q}_s),
\end{equation}
where
\[
\mathbf{q}_i = \mathbf{p}_{j_i}, \quad \text{with } \mathbf{p}_{j_i} \in \mathcal{P}, \quad i = 0, \ldots, s.
\]

The inner mesh \(\mathcal{T}_{\rm in}\) has only an outer boundary, which we denote by
\begin{equation}
\label{eq:GammaInOut}
\Gamma_{\rm in}^{\rm out} := (\mathbf{g}_0, \ldots, \mathbf{g}_r),
\end{equation}
where
\[
\mathbf{g}_i = \mathbf{p}_{k_i}, \quad \text{with } \mathbf{p}_{k_i} \in \mathcal{P}, \quad i = 0, \ldots, r.
\]

The polygonal paths \eqref{eq:GammaOutIn} and \eqref{eq:GammaInOut} are represented as ordered tuples of mesh points, traversed in counterclockwise order. These paths, corresponding to the configuration shown in Fig.~\ref{fig:MeshCoarse}(b), are illustrated in Fig.~\ref{fig:MeshCoarse}(c).

By construction, the paths \(\Gamma_{\rm out}^{\rm in}\) and \(\Gamma_{\rm in}^{\rm out}\) share common points. These points can be readily identified at the combinatorial level by comparing the nodes of the two paths via their global indices in the point set \(\mathcal{P}\).  
As part of the algorithm, we begin by locating the first common point between \(\Gamma_{\rm out}^{\rm in}\) and \(\Gamma_{\rm in}^{\rm out}\). The paths are then synchronized at this point, yielding the initial alignment
\[
\mathbf{q}_0 = \mathbf{g}_0.
\]

If $\mathcal{T}_{\mathrm{in}}$ is well separated from the unstructured boundary layer of the original mesh, we expect an even stronger correlation between the paths \(\Gamma_{\rm out}^{\rm in}\) and \(\Gamma_{\rm in}^{\rm out}\), as illustrated, for example, in Fig.~\ref{fig:MeshCoarse}(c). In particular, we expect that \(\Gamma_{\rm out}^{\rm in}\) contains exactly twice as many points as \(\Gamma_{\rm in}^{\rm out}\), i.e., $s = 2r$, and that the following relation holds between their corresponding nodes:
\begin{equation}
\label{eq:PointsGlued}
\mathbf{q}_{2i} = \mathbf{g}_{i},
\quad \text{for } i = 0, \ldots, r.
\end{equation}

A theoretical guarantee for the existence of a configuration satisfying~\eqref{eq:PointsGlued} is beyond the scope of this work. As part of the algorithm, we check whether condition~\eqref{eq:PointsGlued} is satisfied and terminate the mesh-coarsening procedure if it is not.

As mentioned above, the comparison is performed at the level of global indices in \(\mathcal{P}\), so no round-off ambiguity arises from comparing real-valued positions. Specifically, referring to the path definitions~\eqref{eq:GammaOutIn} and~\eqref{eq:GammaInOut}, we check the condition
\[
j_{2i} = k_i, \quad \text{for } i = 0,  \ldots, r,
\]
which ensures that the correspondence in~\eqref{eq:PointsGlued} holds at the level of global indices in \(\mathcal{P}\).

\subsection{Mesh reconnection}
\label{sec:MeshReconn}

\begin{figure*}[!t]
\centering
\includegraphics[scale=1]{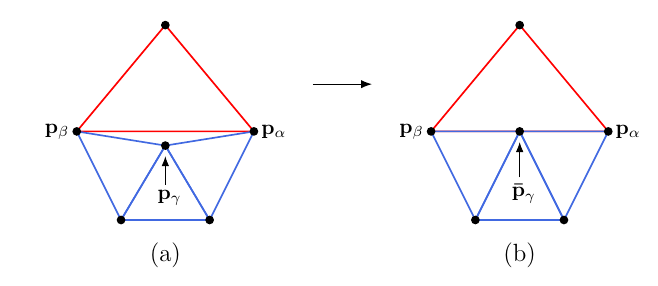}
\caption{\label{fig:DuetContact}
Illustration of the alignment of the contact between the inner and outer meshes. The blue color indicates the elements of the outer mesh, and the red color indicates the elements of the inner mesh. Panel (a) shows a segment of the contact before alignment, and panel (b) shows the same segment after alignment.}
\end{figure*}

Starting from the configuration described in Section~\ref{sec:MeshDecomp}, we merge the inner and outer meshes, \(\mathcal{T}_{\mathrm{in}}\) and \(\mathcal{T}_{\mathrm{out}}\), into a single, simply connected mesh. This is done in two steps.

First, we form an intermediate mesh by taking the union of \(\mathcal{T}_{\mathrm{in}}\) and \(\mathcal{T}_{\mathrm{out}}\):
\[
\mathcal{T}^{*} = \mathcal{T}_{\mathrm{in}} \cup \mathcal{T}_{\mathrm{out}},
\]
with its corresponding set of vertices denoted by \(\mathcal{P}^{*} \subset \mathcal{P}\).

In the second step, we insert void elements into \(\mathcal{T}^{*}\) along the contact region between \(\mathcal{T}_{\mathrm{in}}\) and \(\mathcal{T}_{\mathrm{out}}\), forming nonconforming interfaces that restore the global integrity of the mesh.

To illustrate this procedure, consider the configuration shown in Fig.~\ref{fig:DuetContact}, which resembles that of Fig.~\ref{fig:TriuExamp}. Figure~\ref{fig:DuetContact} shows the contact region between a coarse triangle from the inner mesh and three fine triangles from the outer mesh. The points \(\mathbf{p}_{\alpha}\), \(\mathbf{p}_{\beta}\), and \(\mathbf{p}_{\gamma}\) from \(\mathcal{P}^{*}\) correspond to adjacent nodes located along the paths \(\Gamma_{\rm in}^{\,\rm out}\) and \(\Gamma_{\rm out}^{\rm in}\). Using the notation from Eqs.~\eqref{eq:GammaOutIn}–\eqref{eq:PointsGlued} we consider an index \(i\) such that \(0 \leq i < r\) and assume the following relations:
\begin{align*}
\mathbf{p}_{\alpha} &= \mathbf{q}_{2i+2} = \mathbf{g}_{i+1}, \\
\mathbf{p}_{\beta}  &= \mathbf{q}_{2i}   = \mathbf{g}_{i}, \\
\mathbf{p}_{\gamma} &= \mathbf{q}_{2i+1}.
\end{align*}

We also introduce the constrained node \(\bar{\mathbf{p}}_{\gamma}\), to be used in the construction of the void element:
\begin{equation*}
\bar{\mathbf{p}}_{\gamma} = \tfrac{1}{2} \left( \mathbf{p}_{\alpha} + \mathbf{p}_{\beta} \right).
\end{equation*}

The inclusion of the void element is performed in two steps:
\begin{enumerate}[label=\arabic*)]
\item The position of \(\mathbf{p}_{\gamma}\) is set to \(\bar{\mathbf{p}}_{\gamma}\), without modifying the connectivity matrix of \(\mathcal{T}^{*}\).
\item The row \([\alpha, \beta, \gamma]\), representing the void element, is appended to the connectivity matrix of \(\mathcal{T}^{*}\).
\end{enumerate}
This procedure is repeated for all indices \(i\) in the range \(0 \le i < r\) along the contact region between the inner and outer meshes, finally yielding the new triangle mesh \(\mathcal{T}_{\rm new}\) as the final state of the intermediate mesh \(\mathcal{T}^{*}\).

As discussed in Section~\ref{sec:MeshDef}, the insertion of the nonconforming contacts in this form restores connectivity between mesh elements at the combinatorial level. The resulting mesh \(\mathcal{T}_{\rm new}\) covers a simply connected domain and has a single outer boundary. This mesh can then be used as the input configuration for the next iteration of the coarsening procedure.

\subsection{Algorithm summary}
\label{sec:AlgMeshSumm}

The mesh-coarsening algorithm consisting of the steps described in Sections~\ref{sec:InputConfig},~\ref{sec:MeshDecomp} and \ref{sec:MeshReconn} can be represented as an operator on the set of triangular meshes $\mathcal{M}$. The set consists of tuples $(\mathcal{T}, \mathcal{P})$, where $\mathcal{T}$ is a triangulation over the point set $\mathcal{P}$. The zero element is also included in the set and represents an empty mesh.

The mesh-coarsening algorithm is defined as the operator
\begin{equation}
F^{\rm coarsen}
\left( \alpha_{\rm detach}, \mathbf{c}_0 , n_{\rm shrink}\right) :
\mathcal{M} \to \mathcal{M},
\end{equation}
where $\alpha_{\rm detach} \in \{\text{true}, \text{false}\}$ is the boolean variable that defines if the background supertriangulation is detached from the boundary of the original mesh or not (see Section~\ref{sec:InputBkgSup}), 
$\mathbf{c}_0$ is the anchor point for chosing the starting super-tringle in the reduction step (see Section~\ref{sec:InputIterRed}),
and $n_{\rm shrink}$ is the non-negative integer that defines the number of cleaning steps after the reduction step (see Section~\ref{sec:InputCleaning}).

The operator \(F^{\rm coarsen}\) returns either a new mesh or the zero element when the algorithm terminates. Starting from a given mesh $(\mathcal{T}_0, \mathcal{P}_0)$, the operator can be applied iteratively to produce a sequence of meshes:
\begin{align*}
(\mathcal{T}_1, \mathcal{P}_1) & 
= F^{\rm coarsen}\left[ (\mathcal{T}_0, \mathcal{P}_0) \right],\\
\vdots &\\
(\mathcal{T}_n, \mathcal{P}_n) & 
= F^{\rm coarsen}\left[ (\mathcal{T}_{n-1}, \mathcal{P}_{n-1}) \right],
\end{align*}
with $n$ being the number of iterations.
Each subsequent mesh in the sequence will have a coarse central region generated from the bulk of the previous mesh. As a result, the coarsening procedure, when applied iteratively, generates multiple discretization layers from the input background mesh.

\section{Mesh processing}
\label{sec:MeshProcess}

Meshes generated by the algorithm presented in Section~\ref{sec:MeshConstr} can be further processed to provide additional flexibility in mesh construction.
The unified mesh representation introduced in Section~\ref{sec:MeshDef} enables mesh refinement and coarsening to be performed in a frontal manner.
Refinement and coarsening fronts can be constructed based on the void elements that represent nonconforming element interfaces.

\subsection{Frontal refinement}
\label{sec:FrontRef}

For example, the refinement front can be constructed from triangles that share with the void elements a common edge opposite the hanging node. An example of such a triangle is shown in Fig.~\ref{fig:Fronts}(a). The refinement front can then be further filtered to control the front propagation.

Triangles in the refinement front are refined using the red approach, i.e., by subdividing each triangle into four equal triangles through connecting the midpoints of its edges. The old void element is removed from the mesh, and new void elements are introduced at the newly created hanging nodes to restore the combinatorial integrity of the mesh.

\subsection{Frontal coarsening}
\label{sec:FrontCoarse}

\begin{figure*}
\centering
\includegraphics[scale=0.8]{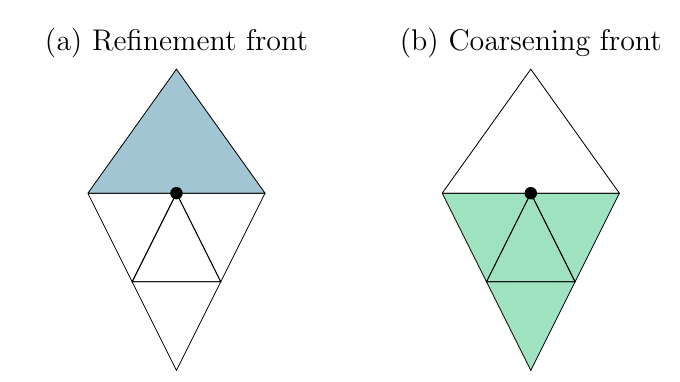}
\caption{\label{fig:Fronts}
Illustration of elements in refinement (a) and coarsening (b) fronts.
The coloured element in panel (a) represents a coarse triangle that is subdivided in the refinement step. The coloured elements in panel (b) represent a group of fine triangles that are replaced by a coarse triangle during the coarsening step.
}
\end{figure*}

The coarsening front can be constructed from supertriangles that contain exactly two subtriangles whose edges are shared with the void elements. An example of such a supertriangle is shown in Fig.~\ref{fig:Fronts}(b). To ensure that coarsening is exactly the inverse of refinement, we also verify that the vertices of the central subtriangle in each supertriangle coincide with the midpoints of the edges. This guarantees that the subtriangles of each supertriangle correspond exactly to those produced by red refinement. Supertriangles that do not satisfy this condition are excluded from the coarsening front. Further filtering can be applied to control front propagation.

The coarsening process proceeds as follows. The subtriangles of each supertriangle are removed from the mesh, and the supertriangles themselves are added to the mesh. The integrity of the mesh is restored by placing void elements at the newly created hanging nodes.

\subsection{Removing defects}
\label{sec:RemovDefects}

Even if a supertriangle in the coarsening front does not consist of subtriangles satisfying the red refinement condition, it can still be coarsened using the procedure described above. Although this process does not possess refinement–coarsening reversibility, it can be applied as a post-processing step after mesh generation. We refer to such supertriangles as \emph{defects}.

The most common defects are supertriangles whose central subtriangle has exactly two hanging nodes as vertices. Such defects may appear at the interface between discretization layers and, in our implementation, are removed in a single coarsening step after the mesh is generated.

\section{Examples}
\label{sec:Examps}

In this section, we present examples of meshes generated by the algorithm described in Section \ref{sec:MeshConstr}. The algorithm has been implemented in the finite-element package \emph{triellipt} \citep{triellipt}.

\subsection{Basic domains}
\label{sec:BasicDomains}

We first demonstrate how the algorithm performs on simple domains.

\begin{figure*}[!h]
\centering
\includegraphics[scale=0.75]{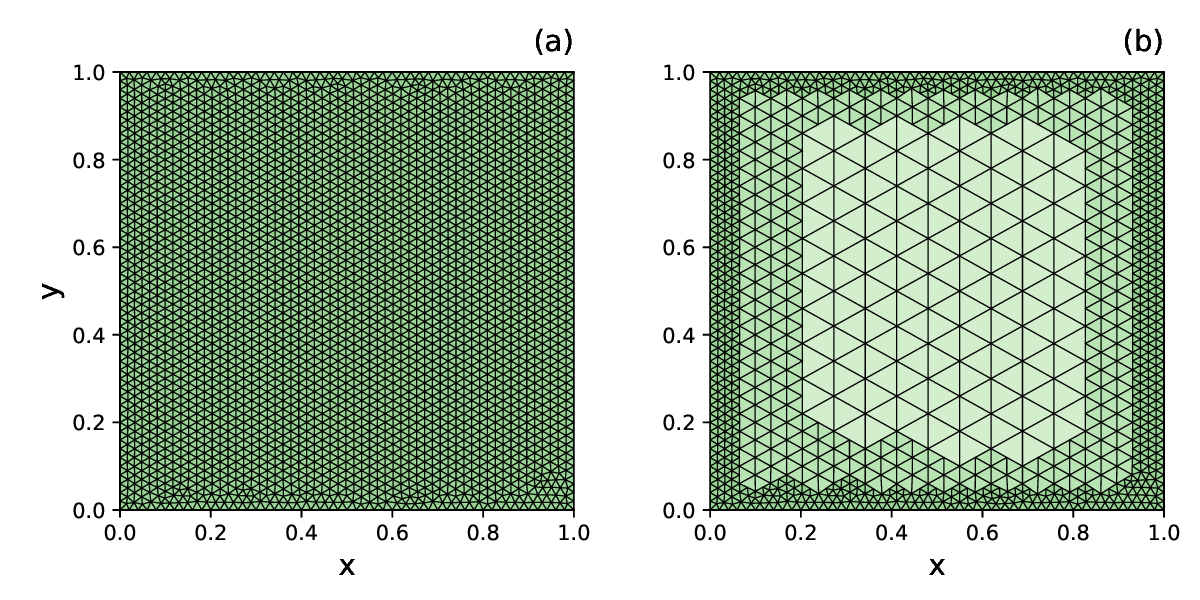}
\caption{\label{fig:RectGmsh}
Example of nonconforming mesh generation in a square domain. Panel (a) shows the input mesh, and panel (b) shows the resulting mesh generated by the proposed framework.}
\end{figure*}

\begin{figure*}[!h]
\centering
\includegraphics[scale=0.75]{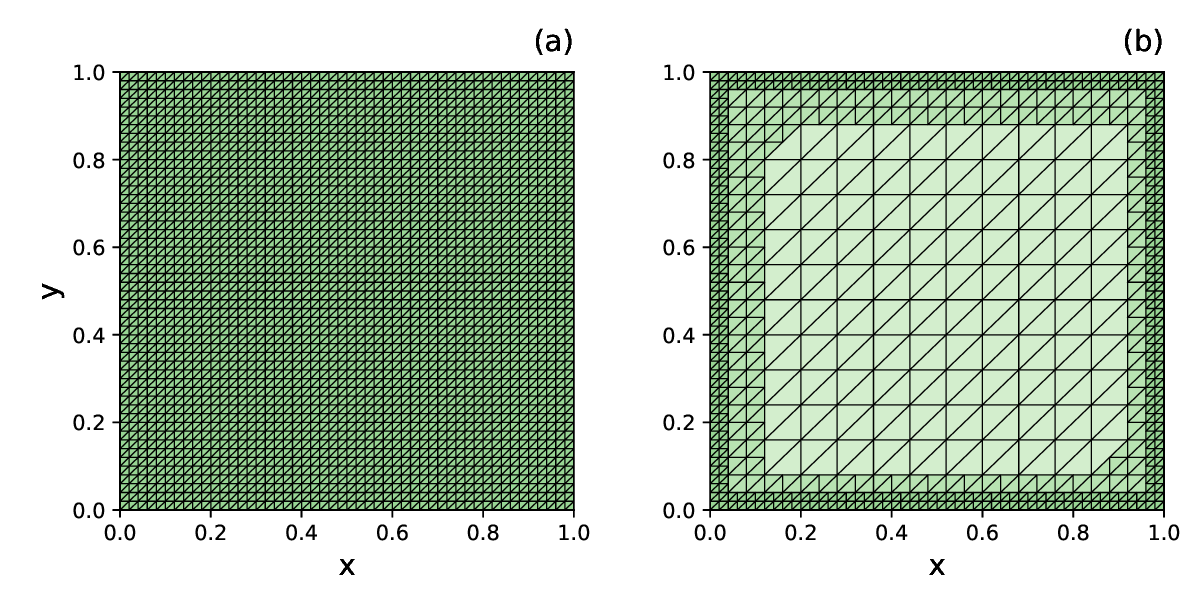}
\caption{\label{fig:RectStruct}
Example of nonconforming mesh generation from a structured triangular mesh.
Panel (a) shows the input mesh, and panel (b) shows the resulting mesh generated by the proposed framework.
}
\end{figure*}

\begin{figure*}[!h]
\centering
\includegraphics[scale=0.75]{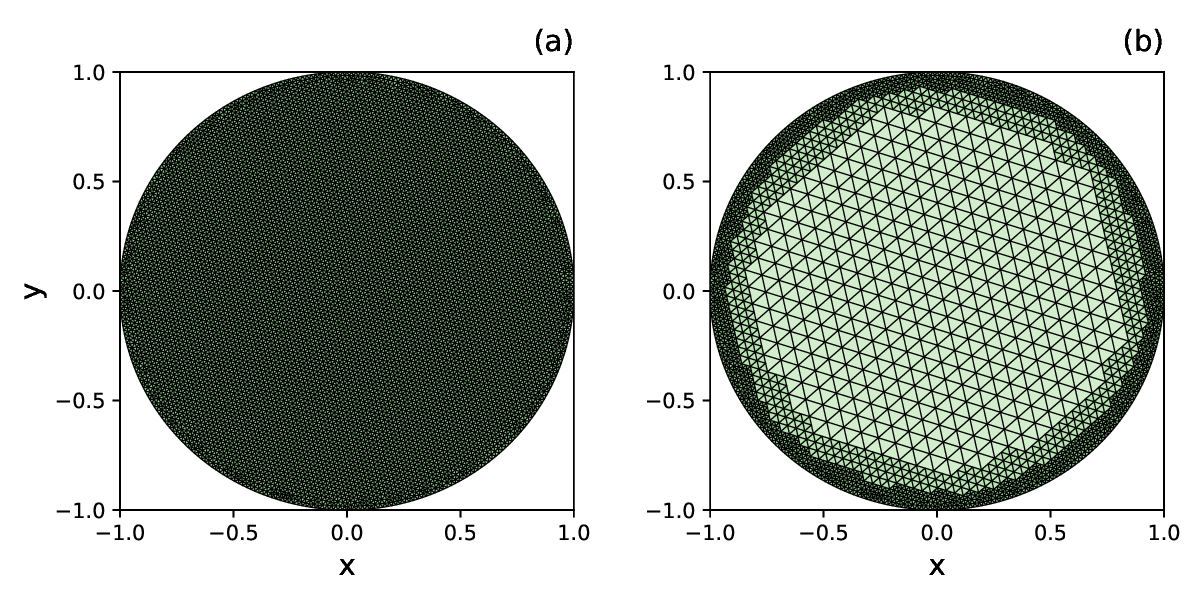}
\caption{\label{fig:CircGmsh}
Example of nonconforming mesh generation in a circular domain.
Panel (a) shows the input mesh, and panel (b) shows the resulting mesh generated by the proposed framework.
}
\end{figure*}

\begin{figure*}[!h]
\centering
\includegraphics[scale=0.70]{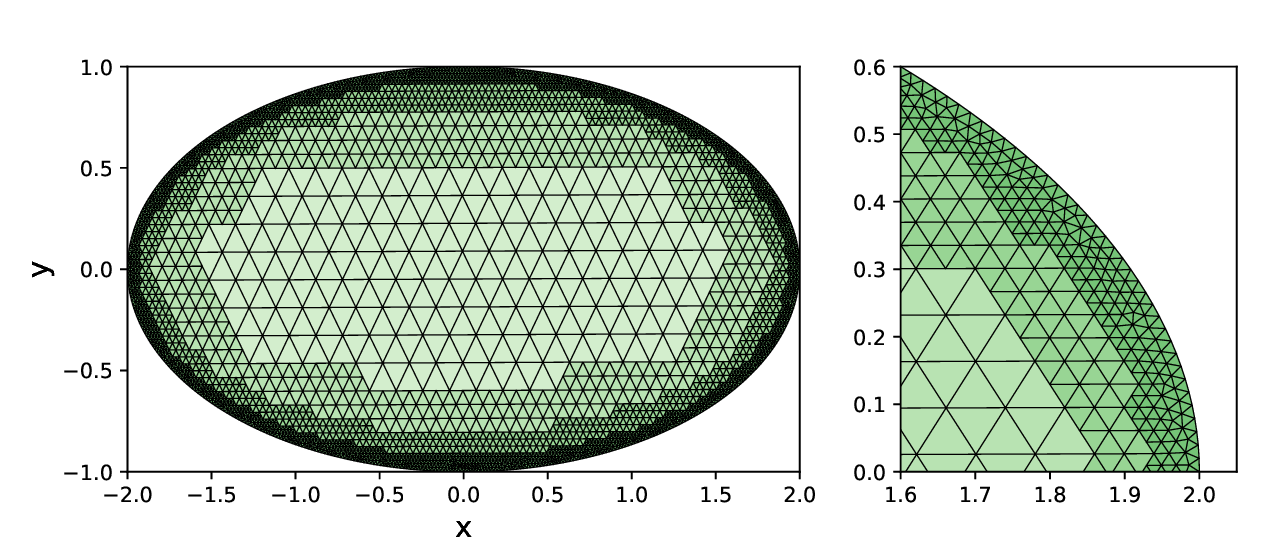}
\caption{\label{fig:EllipGmsh}
Nonconforming mesh with multiple discretization layers generated in an elliptic domain.
}
\end{figure*}

\subsubsection{Unit square}
\label{sec:UnitSquare}

\emph{Unstructured case}.
Figure~\ref{fig:RectGmsh} demonstrates the mesh generation in the unit square. Figure~\ref{fig:RectGmsh}(a) shows the initial background mesh generated with Gmsh~\citep{geuzaine2009gmsh} using the frontal Delaunay algorithm and a uniform mesh size of 0.02. Figure~\ref{fig:RectGmsh}(b) shows a nonconforming mesh generated from the mesh in Fig.~\ref{fig:RectGmsh}(a) by applying the mesh-coarsening algorithm twice.

The first coarsening step was performed with the parameters $\alpha_{\rm detach} = \text{true}$, $\mathbf{c}_0 = (0.5, 0.5)$, and $n_{\rm shrink} = 0$. The second coarsening step was performed with the same parameters, except with $n_{\rm shrink} = 2$.

It can be seen in Fig.~\ref{fig:RectGmsh}(a) that the initial mesh contains an unstructured layer near the boundary of the square and a structured bulk region. Applying the coarsening algorithm to the background mesh separates this structured bulk from the unstructured boundary layer. In this step, changing the parameter $\alpha_{\rm detach}$ can lead to different configurations in subsequent coarsening steps. The parameter $n_{\rm shrink}$ controls the width of the unstructured boundary layer.

The subsequent coarsening steps lead to a coarsening of the structured mesh bulk. In these steps, the parameter $\alpha_{\rm detach}$ does not affect the process, since the mesh-coarsening occurs in the core region, which is already well separated from the boundary. The parameter $n_{\rm shrink}$ controls the separation of the discretization layers in the bulk.

The parameter $\mathbf{c}_0$ determines the starting point of the mesh-coarsening.
Taking into account the algorithm structure discussed above, it is clear that $\mathbf{c}_0$ should be chosen somewhere in the structured bulk to initiate its coarsening. In our implementation, the default value of $\mathbf{c}_0$ is the mean of the centroids of all triangles. This choice is typically sufficient to achieve the expected algorithm performance.

\emph{Structured case}.
The mesh-coarsening algorithm can be applied to certain types of structured triangular meshes. An example of such a case is shown in Fig.~\ref{fig:RectStruct}.

Figure~\ref{fig:RectStruct}(a) shows the initial background mesh, obtained by dividing the cells of a uniform grid along a diagonal. Note that the same diagonal orientation is used in all cells. The grid consists of 51 points in each direction.

Such a mesh possesses a structured supertriangulation that supports the application of the mesh-coarsening algorithm. Figure~\ref{fig:RectGmsh}(b) shows a nonconforming mesh generated from the mesh in Fig.~\ref{fig:RectGmsh}(a) by applying the mesh-coarsening algorithm twice. Both steps were performed with the parameters $\alpha_{\rm detach} = \text{true}$ and $\mathbf{c}_0 = (0.5, 0.5)$. The parameter $n_{\rm shrink}$ was set to 0 in the first step and to 1 in the second step.

The algorithm is capable of generating multiple discretization layers in the considered case. This is not always achievable when other types of structured triangular meshes are used. For example, this situation arises in a cross-triangulated mesh generated by dividing each grid cell along both diagonals. In such a mesh, the supertriangulation consists mostly of disconnected supertriangles, which makes mesh coarsening unfeasible.

\subsubsection{Elliptic domains}
\label{sec:EllipDomains}

We now consider examples of mesh generation in elliptic domains.
We begin with the unit circle case, shown in Fig.~\ref{fig:CircGmsh}.

Figure~\ref{fig:CircGmsh}(a) shows the initial background mesh generated with Gmsh~\citep{geuzaine2009gmsh} using the frontal Delaunay algorithm and a uniform mesh size of 0.02. Figure~\ref{fig:CircGmsh}(b) shows a nonconforming mesh generated from the mesh in Fig.~\ref{fig:CircGmsh}(a) by applying the mesh-coarsening algorithm twice. Both steps were performed with the parameters $\alpha_{\rm detach} = \text{true}$ and $\mathbf{c}_0 = (0, 0)$. The parameter $n_{\rm shrink}$ was set to 1 in the first step and to 2 in the second step.

The qualitative behavior of the algorithm in this case is consistent with that observed in Section~\ref{sec:UnitSquare}. Once again, we observe the separation of the unstructured boundary layer and the structured bulk region, which is then coarsened into multiple discretization layers.

As an additional example, we show in Fig.~\ref{fig:EllipGmsh} a nonconforming mesh generated using the mesh-coarsening algorithm in an elliptic domain with major and minor axes 2 and 1. As in the previous case, the initial background mesh was generated with Gmsh~\citep{geuzaine2009gmsh} using the frontal Delaunay algorithm and a uniform mesh size of~0.02. The nonconforming mesh shown in Fig.~\ref{fig:EllipGmsh} was obtained from the initial mesh by applying the mesh-coarsening algorithm three times. The parameters $\alpha_{\rm detach} = \text{false}$ and $\mathbf{c}_0 = (0,0)$ were used in all three steps. The parameter $n_{\rm shrink}$ was set to~1 in the first step and to~2 in the remaining two steps.

Once again, a clear separation of the unstructured boundary layer and the structured bulk region is visible in Fig.~\ref{fig:EllipGmsh}.

\subsection{Refinement and coarsening}

\begin{figure*}
\centering
\includegraphics[scale=0.75]{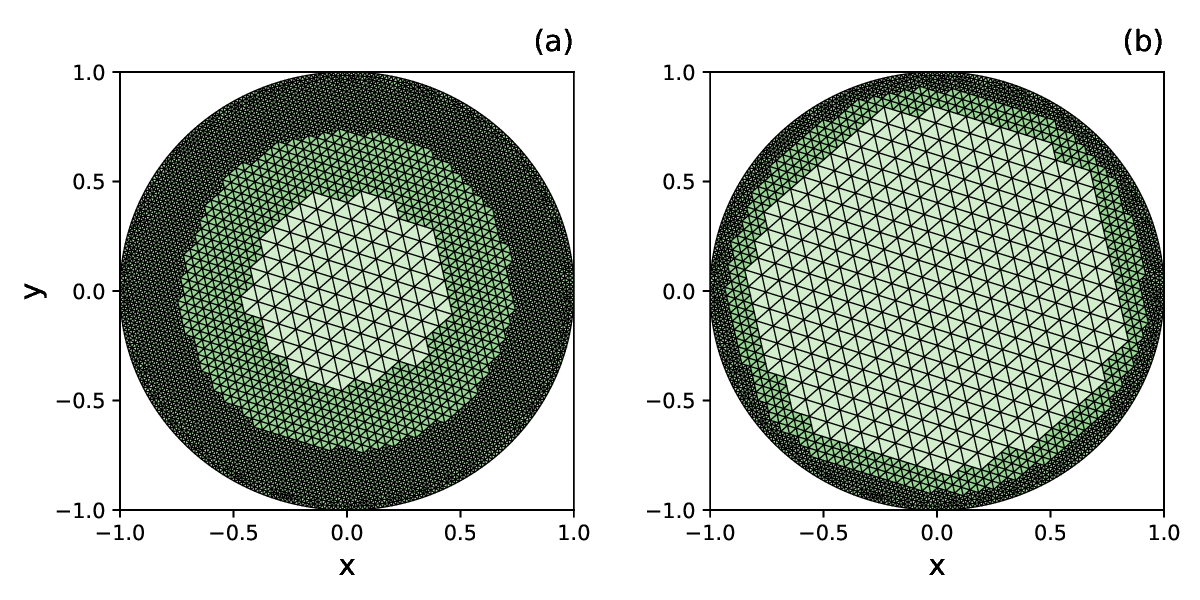}
\caption{\label{fig:CircAMR}
Illustration of frontal refinement and coarsening. Panel (a) shows the mesh generated from the mesh in Fig.~\ref{fig:CircGmsh}(b) after performing ten frontal refinement steps. Panel (b) shows the mesh obtained from the mesh in panel (a) after performing ten frontal coarsening (de-refinement) steps.}
\end{figure*}

As discussed in Section~\ref{sec:MeshProcess}, the nonconforming meshes produced by the mesh-coarsening algorithm can subsequently be refined and de-refined in a frontal manner. To illustrate this, we show in Fig.~\ref{fig:CircAMR} the meshes obtained by frontal refinement and de-refinement starting from the nonconforming mesh shown in Fig.~\ref{fig:CircGmsh}(b).

In this example, we first apply the frontal refinement described in Section~\ref{sec:FrontRef} to the mesh shown in Fig.~\ref{fig:CircGmsh}(b). We perform ten refinement steps, ultimately obtaining the mesh shown in Fig.~\ref{fig:CircAMR}(a). In each refinement step, all triangles in the refinement front are refined.

Starting from the mesh in Fig.~\ref{fig:CircAMR}(a), we apply the de-refinement procedure described in Section~\ref{sec:FrontCoarse}. We perform ten derefinement steps ultimately obtaining the mesh shown in Fig.~\ref{fig:CircAMR}(b). In each coarsening step, all triangles in the coarsening front are coarsened. 

As expected, the mesh in Fig.~\ref{fig:CircAMR}(b) coincides with the initial mesh shown in Fig.~\ref{fig:CircGmsh}(b). This is because we applied the same number of refinement and de-refinement steps. This confirms the reversibility of the refinement and coarsening operations described in Section~\ref{sec:MeshProcess}.

The example shown in Fig.~\ref{fig:CircAMR} demonstrates that the nonconforming meshes considered in this work provide sufficient flexibility for adaptive mesh refinement. Frontal refinement and de-refinement offer a convenient alternative to conventional adaptive mesh-refinement techniques, such as those based on binary-tree or octree data structures.

\begin{figure*}
\centering
\includegraphics[scale=0.75]{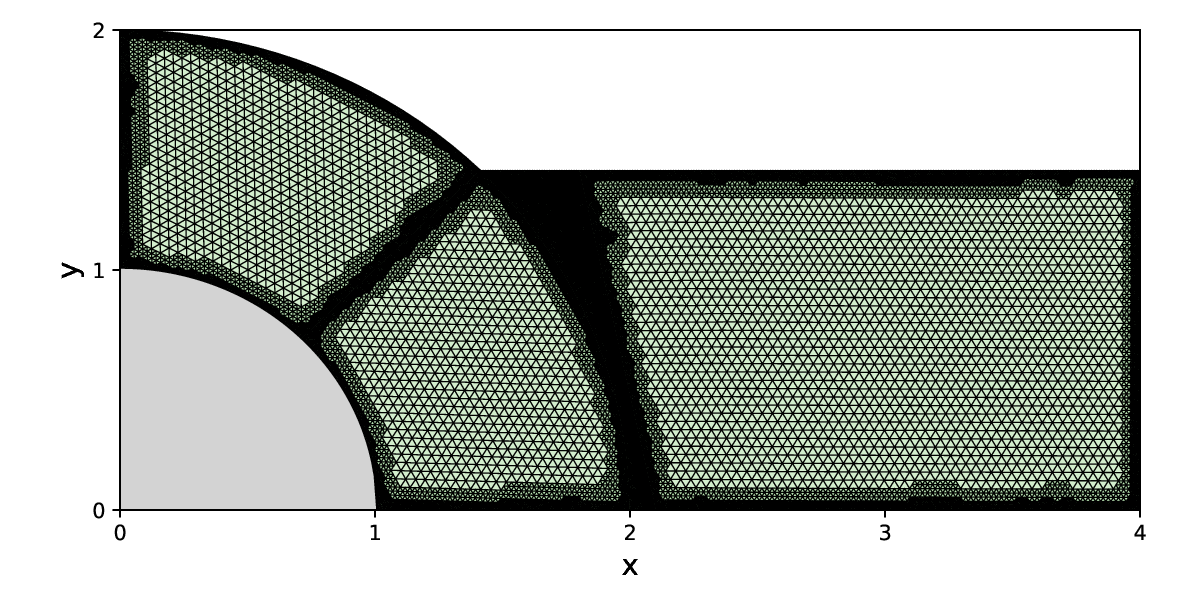}
\caption{\label{fig:Composite}
Example of a composite multidomain nonconforming mesh
with multiple discretization layers.}
\end{figure*}

\subsection{Composite mesh}

In Section~\ref{sec:BasicDomains}, we applied the mesh-coarsening algorithm to simply shaped domains. Applying this algorithm to more complex geometries is not straightforward, as the main difficulty lies in generating a background mesh of sufficiently high quality.

Meshes generated with the frontal Delaunay algorithm can exhibit well-structured bulk regions even in complex settings. These regions can be separated by the mesh-coarsening algorithm through adjustments of the seed parameter~\(\mathbf{c}_0\). However, the overall mesh quality in the bulk region is lower than for simply shaped domains, which reduces algorithm performance and necessitates more careful parameter adjustment.

For complex geometries, the proposed mesh generation framework can be used to construct composite multidomain meshes. In this approach, the domain is first decomposed into subdomains of simpler geometry, and each subdomain is subsequently discretized using a mesh with multiple discretization layers. An example of such a mesh is shown in Fig.~\ref{fig:Composite}.

This type of geometry arises, for example, in studies of atmospheric-pressure gas discharges in the so-called rod-to-plane configuration \citep{takahashi1994, ferreira2019, ren2022}. From a physical point of view, such discharges exhibit a complex spatio-temporal behavior, where thin ionization waves can originate at the rod tip and propagate through the domain. Effective simulation strategies for these discharges typically involve the use of adaptive mesh refinement, which is mostly applied on structured grids for open-space problems \citep{teunissen2017, semenov2022}.

The transition to conservative adaptive mesh refinement on unstructured meshes requires a discretization of the domain that is suitable for applying red refinement and de-refinement procedures. Meshes like the one shown in Fig.~\ref{fig:Composite} are a possible choice for such a discretization. In such a mesh, the coarse bulk regions within the subdomains can be refined and de-refined in a controlled and well-defined manner. Moreover, the interpolation operators can be constructed independently for each subdomain. Such meshes may therefore be useful in constructing effective simulation strategies for multiscale problems.

\section{Conclusion}
\label{sec:Concls}

We present a framework for generating nonconforming triangular meshes with multiple discretization layers. The framework takes as input background conforming meshes generated by the frontal Delaunay algorithm with a uniform element size. Owing to the characteristic structure of such meshes, their bulk can be coarsened in a repetitive manner, resulting in nonconforming meshes with a hierarchy of discretization layers.

The element size in the mesh increases from the boundary toward the bulk. The finest discretization layer is responsible for resolving the boundary of the discretized domain. The remaining discretization layers, with increasing element size, ensure a transition from the unstructured boundary layer to the structured bulk region of the mesh.

The performance of the proposed algorithm is influenced by the quality of the input background mesh. The expected mesh structure is generally achieved for domains with relatively simple geometry. For more complex geometries, the proposed framework can be employed to construct composite multidomain meshes, where a hierarchy of discretization layers is generated within each subdomain.

Meshes generated using the proposed framework are well suited for adaptive mesh refinement techniques. Mesh refinement and coarsening can be performed locally or in a frontal manner. In the latter case, the mesh is processed by refining or coarsening the elements at the interfaces between the discretization layers.

The proposed mesh generation framework is straightforward to implement in combination with modern mesh generators such as Gmsh \cite{geuzaine2009gmsh}. It is currently implemented in the finite-element package {\it triellipt} \citep{triellipt}. The framework can also be readily integrated with other finite-element packages that support nonconforming triangular meshes.

\section*{Acknowledgment}
This work was funded by the Deutsche Forschungsgemeinschaft (DFG, German Research Foundation) – project number 515939493.



\end{document}